\documentclass[aap,MSNbibl,citesort,dvips]{arximspdf}

\doi{10.1214/11-AAP798}
\volume{22}
\issue{4}
\pubyear{2012}
\firstpage{1450}
\lastpage{1464}

\makeatletter
\newtheorem{theorem}{Theorem}
\newtheorem{lemma}[theorem]{Lemma}
\newtheorem{corollary}[theorem]{Corollary}

\newproclaim{remark}[theorem]{Remark}

\newcommand{\eps}{\varepsilon}
\renewcommand{\Pr}{{\mathbb P}}
\newcommand{\E}{{\mathbb E}}

\newcommand{\cF}{\mathcal{F}}
\newcommand{\cC}{\mathcal{C}}

\newcommand{\floor}[1]{\lfloor #1 \rfloor}
\newcommand{\ceil}[1]{\lceil #1 \rceil}

\newcommand{\bv}{{\underline v}}
\newcommand{\pto}{\stackrel{\mathrm{p}}{\to}}
\newcommand{\cR}{{\mathcal R}}
\newcommand{\tc}{t_{\mathrm{c}}}
\newcommand{\cL}{{\mathcal L}}
\newcommand{\cX}{{\mathcal X}}
\newcommand{\cG}{{\mathcal G}}

\makeatother

\begin{document}
\begin{frontmatter}

\title{Achlioptas process phase transitions are continuous}
\runtitle{Achlioptas process phase transitions are continuous}

\begin{aug}
\author[A]{\fnms{Oliver} \snm{Riordan}\corref{}\ead[label=e1]{riordan@maths.ox.ac.uk}}
\and
\author[A]{\fnms{Lutz} \snm{Warnke}\ead[label=e2]{warnke@maths.ox.ac.uk}}
\runauthor{O. Riordan and L. Warnke}
\affiliation{University of Oxford}
\address[A]{Mathematical Institute\\
University of Oxford\\
24--29 St Giles'\\
Oxford OX1 3LB\\
United Kingdom\\
\printead{e1}\\
\phantom{E-mail: }\printead*{e2}} 
\end{aug}

\received{\smonth{3} \syear{2011}}
\revised{\smonth{8} \syear{2011}}

%
\begin{abstract}
It is widely believed that certain simple modifications of the random
graph process lead to discontinuous phase transitions. In particular,
starting with the empty graph on $n$ vertices, suppose that at each
step two pairs of vertices are chosen uniformly at random, but only one
pair is joined, namely, one minimizing the product of the sizes of the
components to be joined. Making explicit an earlier belief of
Achlioptas and others, in 2009, Achlioptas, D'Souza and Spencer
[\textit{Science} \textbf{323} (2009) 1453--1455] conjectured that
there exists a $\delta>0$ (in fact, $\delta\ge1/2$) such that with
high probability the order of the largest component ``jumps'' from
$o(n)$ to at least $\delta n$ in $o(n)$ steps of the process, a
phenomenon known as ``explosive percolation.''

We give a simple proof that this is not the case. Our result applies to
all ``Achlioptas processes,'' and more generally to any process where a
fixed number of independent random vertices are chosen at each step,
and (at least) one edge between these vertices is added to the current
graph, according to any (online) rule.

We also prove the existence and continuity of the limit of the rescaled
size of the giant component in a class of such processes, settling a
number of conjectures. Intriguing questions remain, however, especially
for the product rule described above.
\end{abstract}

%

\vspace*{-2pt}
\end{frontmatter}

\section{Introduction and results}\label{secintro}

At a Fields Institute workshop in 2000, Dimitris Achlioptas suggested a
class of variants of the classical random graph process, defining a
random sequence $(G(m))_{m\ge0}$ of graphs on a fixed vertex set of
size~$n$, usually explained in terms of the actions of a hypothetical
purposeful agent: start at step $0$ with the empty graph. At step $m$,
two potential edges $e_1$ and $e_2$ are chosen independently and
uniformly at random from all ${n\choose2}$ possible edges [or from
those edges not present in $G(m-1)$]. The agent must select one of
these edges, setting $G(m)=G(m-1)\cup\{e\}$ for $e=e_1$ or $e_2$. Any
possible strategy, or ``rule,'' for the agent gives rise to a~random
graph process. Such processes are known as ``Achlioptas
processes.''\looseness=-1\vadjust{\goodbreak}

If the agent always chooses the first edge, then (ignoring the minor
effect of repeated edges) this is, of course, the classical random
graph process, studied implicitly by Erd\H os and R\'enyi and
formalized by Bollob\'as. In this case, as is well known,
there is a phase transition around $m=n/2$.
More precisely, writing $L_1(G)$ for the number of vertices
in the (a, if there is a~tie) largest component of a graph~$G$,
Erd\H os and R\'enyi~\cite{ERgiant} showed that
there is a function $\rho=\rho^{\mathrm{ER}}\dvtx [0,\infty)\to[0,1)$ such
that for any fixed $t\ge0$, whenever $m=m(n)$ satisfies $m/n\to t$ as
$n\to\infty$,
then $L_1(G(m))/n \pto\rho(t)$, where $\pto$ denotes
convergence in probability. Moreover, $\rho(t)=0$ for $t\le1/2$,
$\rho(t)>0$ for $t>1/2$ and $\rho(t)$ (the solution
to a simple equation) is continuous at $t=1/2$ with
right-derivative $4$ at this point.

Achlioptas originally asked whether the agent could shift the critical
point of this phase transition by following an appropriate
edge-selection rule. One
natural rule to try is the ``product rule'': of the given potential
edges, pick the one minimizing the product of the sizes of the
components of its endvertices. This rule was suggested by Bollob\'as as
the most likely to delay the critical point.

Bohman and Frieze~\cite{BF2001} quickly showed, using a much simpler
rule, that the
transition could indeed be shifted, but more complicated rules such as
the product rule remained resistant to analysis. By 2004 at the
latest (see~\cite{SpencerWormald2007}),
extensive simulations of D'Souza and others strongly suggested that
the product rule in particular shows much
more interesting behavior than simply a slightly shifted critical
point; it exhibits a phenomenon known as ``explosive percolation.''

As usual, we say that an event $E$ (formally a sequence of events $E_n$)
holds \textit{with high probability} (\textit{whp}) if $\Pr(E)\to1$
as $n\to\infty$.
Explosive percolation is said to occur if there is a critical $\tc$ and
a positive
$\delta$ such that for any fixed $\eps>0$,
whp $L_1$ jumps from $o(n)$ to at least $\delta n$ in fewer than $\eps n$
steps around $m=\tc n$.
Recently, Achlioptas, D'Souza and Spencer~\cite{AchlioptasDSouzaSpencer2009}
presented ``conclusive numerical evidence'' for the conjecture that
the product rule exhibits explosive percolation,
suggesting indeed that the largest component
grows from size at most $\sqrt{n}$ to size at least $n/2$
in at most $2n^{2/3}$ steps. Bohman~\cite{Bohman2009} describes
this explosive percolation conjecture as an important and intriguing
mathematical
question.

Our main result disproves this conjecture. The result applies to all
Achlioptas processes as defined at the start of the section (including
the product rule) and, in fact, to a more general class of processes
($\ell$-vertex rules) defined in Section~\ref{secproof}. A~form of this
result first appeared in~\cite{Science}, with more restrictive
assumptions, and without full technical details.
%
%
\begin{theorem}\label{th1}
Let $\cR$ be an $\ell$-vertex rule for some $\ell\ge2$.
For each $n$, let $(G(m))_{m\ge0}=(G_n^{\cR}(m))_{m\ge0}$ be the random
sequence of graphs on $\{1,2,\ldots,n\}$ associated to $\cR$.
Given any functions $h_L(n)$ and $h_m(n)$ that are $o(n)$,
and any constant $\delta>0$,
the probability that there exist $m_1$ and $m_2$
with\break $L_1(G(m_1))\le h_L(n)$, $L_1(G(m_2))\ge\delta n$ and $m_2\le m_1+h_m(n)$
tends to $0$ as $n\to\infty$.
\end{theorem}

Let $N_k(G)$ denote the number of vertices of a graph $G$ in components
with $k$ vertices, so $N_k(G)$ is $k$ times the number of $k$-vertex components.
Similarly, $N_{\le k}(G)$ and $N_{\ge k}(G)$ denote the number of
vertices in components
with at most (at least) $k$ vertices.
Having a rule $\cR$ in mind, and suppressing
the dependence on~$n$, we write $N_k(m)$ for the
random quantity $N_k(G(m))$, and similarly $L_1(m)$ for $L_1(G(m))$.

Under a mild additional condition (which holds for all Achlioptas processes),
a slight modification of the proof of Theorem~\ref{th1} shows, roughly
speaking, that the giant component is unique. In fact, we obtain much more;
whp there is no time at which there are ``many'' vertices in ``large''
components but not in the single largest component. For the precise definition
of a ``merging'' rule see Section~\ref{secedge}; any Achlioptas process
is merging.
%
%
\begin{theorem}\label{thgu}
Let $\cR$ be a merging $\ell$-vertex rule for some $\ell\ge2$.
For each $n$, let $(G(m))_{m\ge0}=(G_n^{\cR}(m))_{m\ge0}$ be the random
sequence of graphs on $\{1,2,\ldots,n\}$ associated to $\cR$.
For each $\eps>0$ there is a $K=K(\eps,\ell)$ such that
\[
\Pr\bigl(\forall m\dvtx N_{\ge K}(m) < L_1(m)+\eps n\bigr) \to1
\]
as $n\to\infty$.
\end{theorem}

With $\ell$ fixed, our proof gives a value for $K$ of the form $\exp
(\exp(c \eps^{-(\ell-1)}))$
for some positive $c=c(\ell)$. Furthermore, we can allow $\eps$ to
depend on $n$,
as long as $\eps=\eps(n)\ge d/(\log\log n)^{1/(\ell-1)}$, where
$d=d(\ell)>0$.

For the classical random graph process it is well known that at any
fixed time,
whp there will be at most one ``giant'' component. Indeed, the maximum
size of
the second largest component throughout
the evolution of the process is whp $o(n)$; this can be read out
of the original results of Erd\H os and R\'enyi~\cite{ERgiant} or
(more easily) the more precise results of Bollob\'as~\cite{BB-evo}.
Spencer's ``no two giants'' conjecture (personal communication) states
that this should also hold for Achlioptas processes.
Theorem~\ref{thgu} proves this conjecture for the larger class of
merging $\ell$-vertex rules; indeed, it readily implies that, with high
probability, the second largest component has size at most
$\max\{K,\eps n\} = \eps n$. Allowing $\eps$ to vary with $n$ as noted
above, the bound we obtain is of the form $d(\ell) n/(\log\log
n)^{1/(\ell-1)}$.

Before turning to the proofs of Theorems~\ref{th1} and~\ref{thgu},
let us discuss some related questions of convergence.

We say that the rule $\cR$ is \textit{locally convergent} if
there exist functions $\rho_k=\rho_k^{\cR}\dvtx[0,\infty)\to[0,1]$
such that,
for each fixed $k\ge1$ and $t\ge0$, we have
%
%
\begin{equation}\label{ktconv}
\frac{N_k(\floor{tn})}{n} \pto\rho_k(t)
\end{equation}
as $n\to\infty$. The rule $\cR$ is \textit{globally convergent} if
there exists an increasing
function $\rho=\rho^{\cR}\dvtx[0,\infty)\to[0,1]$ such that for any $t$
at which $\rho$ is continuous we have
\[
\frac{L_1(\floor{tn})}{n} \pto\rho(t)
\]
as $n\to\infty$.

Theorem~\ref{th1} clearly implies that if a rule $\cR$ is globally
convergent, then the limiting function $\rho$ is continuous at
the critical point $\tc=\inf\{t\dvtx\rho(t)>0\}$. Using Theorem~\ref{thgu},
it is not hard to establish continuity elsewhere for merging rules; see
Theorem~\ref{thcont} and Corollary~\ref{corcont} in Section~\ref{secedge}.
Unfortunately, we cannot show that the product rule \textit{is} globally
convergent. However, as we shall see in Section~\ref{secconv},
Theorem~\ref{thgu} implies the following result.
%
%
\begin{theorem}\label{thlg}
Let $\cR$ be a merging $\ell$-vertex rule for some $\ell\ge2$.
If $\cR$ is locally convergent, then $\cR$ is globally convergent,
and the limiting function~$\rho^{\cR}$ is continuous and satisfies
$\rho^{\cR}(t)=1-\sum_{k\ge1}\rho^{\cR}_k(t)$.
\end{theorem}

The conditional result above is, of course, rather unsatisfactory.
However, for many Achlioptas processes, local convergence is well
known; global convergence had not previously been established for any
nontrivial rule. In particular, Theorem~\ref{thlg} settles two
conjectures of Spencer and Wormald~\cite{SpencerWormald2007} concerning
so-called ``bounded size Achlioptas processes'' (see
Section~\ref{size}).

Recently, in a paper in the physics literature,
da Costa, Dorogovtsev, Goltsev and Mendes~\cite{dCDGM}
announced a version of Theorem 1. However, their actual analysis
concerned only one specific rule (not the product rule, though they
claim that ``clearly'' the product rule is less likely to have a
discontinuous transition). More importantly, even the ``analytic'' part
of it is heuristic, and of a type that seems to us very hard (if at
all possible) to make precise. Crucially, the starting point for their
analysis is not only to assume convergence, but also to assume that
the phase transition is continuous! From this, and some further
assumptions, by solving approximations to certain equations they
deduce certain ``self-consistent behavior,'' which apparently justifies
the assumption of continuity. The argument (which is considerably more
involved than the simple proof presented here) is certainly very
interesting, and the conclusion is (as we now know) correct, but it
seems to be very far from a mathematical proof.

In the next section we prove Theorem~\ref{th1}.
In Section~\ref{secedge}, restricting
the class of rules slightly, we prove Theorem~\ref{thgu}
and deduce that jumps in $L_1$ are also impossible
after a giant component first emerges.
Next, in Section~\ref{secconv},
we prove Theorem~\ref{thlg}.
Finally, in Section~\ref{size} we consider more
restrictive rules such as bounded size rules,
and discuss the relationship of our results to earlier work.

\section{\texorpdfstring{Definitions and proof of Theorem \protect\ref{th1}}{Definitions and proof of Theorem 1}}\label{secproof}

Throughout, we fix an integer \mbox{$\ell\ge2$}. For each $n$, let
$(\bv_1,\bv_2,\ldots)$ be an i.i.d. sequence where each
$\bv_m$ is a sequence $(v_{m,1},\ldots,v_{m,\ell})$ of $\ell$ vertices
from $[n]=\{1,2,\ldots,n\}$ chosen independently and uniformly at
random. Suppressing the dependence on $n$,
informally, an \textit{$\ell$-vertex rule} is a random sequence
$(G(m))_{m\ge0}$ of graphs on~$[n]$ satisfying (i) $G(0)$ is the empty graph,
(ii) for $m\ge1$ $G(m)$ is formed
from $G(m-1)$ by adding a (possibly empty) set $E_m$
of edges, with all edges in~$E_m$ between vertices in $\bv_m$
and (iii) if all $\ell$ vertices in $\bv_m$ are in distinct components
of $G(m-1)$, then $E_m\ne\varnothing$. The set $E_m$ may be chosen
according to any deterministic or random online rule.

Formally, we assume the existence of a filtration $\cF_0\subseteq\cF
_1\subseteq\cdots$
such that~$\bv_m$ is $\cF_m$-measurable and independent
of $\cF_{m-1}$, and require $E_m$ [and hence,~$G(m)$] to
be $\cF_m$-measurable.

In other words, the agent is presented with the random list (set) $\bv_m$
of vertices, and, unless two or more are already in the same component,
must add one or more edges between them, according to any deterministic
or random rule that depends only on the history. In the original examples
of Achlioptas, the rule always adds either the edge $\{v_{m,1},v_{m,2}\}$
or the edge $\{v_{m,3},v_{m,4}\}$. Note that (for now) no connection
between the algorithms
used for different $n$ (or indeed at different steps $m$) is assumed.

The arguments that follow are robust to small changes
in the definition, since
they can be written to rely only on deterministic properties
of $(G(m))$, plus bounds on the probabilities
of certain events at each step.
The latter always have $\Theta(1)$ elbow room.
It follows that we may weaken
the conditions on $(\bv_m)$; it suffices if, for $m=O(n)$, say, the
conditional distribution
of $\bv_m$ given the history (i.e., given $\cF_{m-1}$)
is close to [at total variation distance~$o(1)$ from,
as $n\to\infty$] that described above. This covers variations
such as picking an $\ell$-tuple of \textit{distinct} vertices, or picking
(the ends of) $\ell/2$ randomly selected (distinct) edges not already
present in $G(m-1)$.

The proof of Theorem~\ref{th1} is based on two observations,
which we first present in heuristic form.

Observation 1: If at some time $t$ (i.e., when $m\sim tn$) there are
$\alpha n$ vertices in components of order at least $k$, then within
time $\gamma=O(1/(\alpha^{\ell-1}k))$ a~component of order at least
$\alpha n/\ell^2=\beta n$ will emerge. Indeed, fix a set
$W$ with $|W|\ge\alpha n$ consisting of components of order at least
$k$. At every subsequent step we have probability at least
$\alpha^\ell$ of choosing only vertices in $W$, and if no component
has order more than $\beta n$, it is likely that all these vertices
are in different components, so the rule is \textit{forced} to join two
components meeting~$W$. This cannot happen more than $|W|/k$
times.\vadjust{\goodbreak}

(A form of Observation 1 appears in a paper
of Friedman and Landsberg~\cite{FriedmanLandsberg2009} as a key
part of a heuristic argument \textit{for} explosive percolation.
It is not quite stated correctly, although this does not seem
to be why the heuristic fails.)

Observation 2: Components of order $k$
have a half-life that may be bounded in terms of $k$;
in an individual step, such a component disappears (by joining
another component) with probability at most $k\ell/n$.
Assuming (which we shall not assume in the actual proof)
that the rule $\cR$ is locally convergent, it follows easily
that for all $t_1$, $t_2$ and $k$ we have
$\rho_k(t_1+t_2)\ge\rho_k(t_1)e^{-k\ell t_2}$.

We place vertices into ``bins'' corresponding to component sizes
between~$2^j$ and $2^{j+1}-1$,
writing $\sigma_j(t)$ for $\sum_{2^j\le k<2^{j+1}} \rho_k(t)$.
Let $\alpha>0$ be constant and suppose that $\sigma_j(t)\ge\alpha$ for
some $t<\tc$.
Writing $k=2^j$,
by Observation 1 we have $\tc-t=O(1/k)$, with the implicit constant
depending on $\alpha$, since the $\ge\alpha n$
vertices in components of size at least $k$ will quickly form
a giant component. Using Observation 2, it follows that $\sigma_j(\tc
)\ge g(\alpha)>0$,
for some (explicit but irrelevant) function $g(\alpha)$.

Let $\sigma_j=\sup_{t\le\tc}\sigma_j(t)$.
If $\sigma_j>\alpha$, then $\sigma_j(\tc)\ge g(\alpha)$.
Counting vertices, we have $\sum_j \sigma_j(\tc)\le1$.
Hence, for each $\alpha>0$, only a finite number of $\sigma_j$ can
exceed~$\alpha$.
Thus $\sigma_j\to0$ as $j\to\infty$.
It follows that for any constant $B\ge2$ and any $k=k(n)\to\infty$,
at no $t=t(n)<\tc$ can there be $\Theta(n)$ vertices
in components of size between $k$ and $Bk$.

Using Observation 1, it is easy to deduce that there cannot
be a discontinuous transition. Indeed, if $\lim_{t\to\tc^{ +}}\rho(t)
\ge\delta>0$,
then for any $k$, at time $t_k=\tc-\delta/(\ell^2 k)$, there must be
at least
$\delta n/2$ vertices in components of order at least $k$, so $\rho
_{\ge k}(t_k)\ge\delta/2$,
where $\rho_{\ge k}=1-\sum_{k'<k}\rho_{k'}$.
For any constant $B\ge2$, if $k$ is large it follows that $\rho_{\ge
Bk}(t_k)\ge\delta/3$.
Taking $B$ large enough,
Observation 1 then implies that $\tc-t_k$ is much smaller than $\delta
/(\ell^2k)$.

We now make the above argument precise, without assuming convergence. This
introduces some minor additional complications, but they are easily handled.
We start with two lemmas corresponding to the two observations above.
%
%
\begin{lemma}
\label{lemcreateGC}
Given $0 < \alpha\le1$, let $\cC(\alpha)$ denote the event that for
all $0 \le m \le n^2$
and $1 \le k \le\frac{\alpha}{16} \frac{n}{\log n}$ the following
holds: $N_{\ge k}(m) \ge\alpha n$ implies $L_1(m+\Delta) > \frac
{\alpha
}{\ell^2} n$ for $\Delta= \ceil{\frac{4}{\alpha^{\ell-1}}\frac{n}{k}}$.
Then $\Pr(\cC(\alpha))\ge1-n^{-1}$.
\end{lemma}
\begin{pf}
It suffices to consider fixed $m$ and $k$ and show that, conditional
on $\cF_m$, if $G(m)$ satisfies $N_{\ge k}(m) \ge\alpha n$,
then we have $L_1(m+\Delta) > \frac{\alpha}{\ell^2} n$ with probability
at least $1-n^{-4}$.

Condition on $\cF_m$. Let $W$ be the union of all components with size
at least~$k$ in $G(m)$,
set $\tilde{\alpha} = |W|/n \ge\alpha$ and let $\beta= \tilde
{\alpha
}/\ell^2$.
We now consider the next~$\Delta$ steps.

We say that a step is \textit{good} if (a) all $\ell$ randomly chosen
vertices are in~$W$ and (b) all these vertices are in different components.
Let $X_j$ denote the indicator function of the event that step $m+j$ is good.
Set $X = \sum_{1 \le j \le\Delta} X_{j}$ and $Y = \sum_{1 \le j \le
\Delta} Y_{j}$, where
\[
Y_{j}   =   \cases{
X_{j}, &\quad if $L_1(m+j-1) \le\beta n$,\cr
1, &\quad otherwise.}
\]
Clearly, in each step (a) holds with probability $\tilde{\alpha}^\ell$.
Furthermore, whenever $L_1(m+j-1) \le\beta n$ holds, in step $m+j$ the
probability
that (a) holds and (b) fails is at most ${\ell\choose2}\tilde{\alpha
}^{\ell-1} \beta< \tilde{\alpha}^\ell/2$
(there must be $v_a$ and $v_b$ with $1\le a < b \le\ell$ such that
$v_b$ lies
in the same component as $v_a$; all $v_c$ must also be in $W$)
and so in this case step $m+j$ is good with probability at least
$\tilde{\alpha}^\ell/2$.
Since, otherwise, $Y_{j}=1$ by definition, we deduce that $Y$
stochastically dominates a binomial random variable with mean $\Delta
\tilde{\alpha}^\ell/2 \ge2 \tilde{\alpha} n /k$.
Standard Chernoff bounds now imply that $\Pr(Y \le\tilde{\alpha} n
/k)\le e^{-\tilde{\alpha} n /(4k)} \le e^{-\alpha n /(4k)} \le n^{-4}$.

Assume that $L_1(m+\Delta) \le\beta n$. Then by monotonicity
$L_1(m+j-1) \le\beta n$ for every $1 \le j \le\Delta$, so $X = Y$.
Note that $W$ contains at most $|W|/k = \tilde{\alpha} n/k$ components
in $G(m)$.
Since every good step joins two components meeting~$W$
[at least one such edge must be added since by (a) all endpoints are in
$W$ and by (b) all endpoints are in distinct components]
we deduce that $Y \le\tilde{\alpha} n/k$.
Hence, $\Pr(L_1(m+\Delta) \le\beta n) \le\Pr(Y \le\tilde{\alpha
} n
/k) \le n^{-4}$, as required.
\end{pf}

Applying Lemma~\ref{lemcreateGC} with $m=0$, $k=1$ and $\alpha=1$, we
readily deduce that whp a giant component exists after at most $4n$ steps.
In fact, it is easy to see that for any $\eps>0$,
whp there is a giant component after at most $(1+\eps)n$ steps (see
the proof of Lemma~\ref{cGC2}).
%
%
\begin{lemma}
\label{lemBinSize}
Fix $0 < \alpha\le1$, $D>0$ and an integer $B\ge2$.
Define $M_k^B(m) = N_{\ge k}(m)-N_{\ge Bk}(m)$.
Let $\cL(\alpha,B,D)$ denote the event that for all $0 \le m \le n^2$
and $1 \le k \le\min\{\frac{\alpha^2e^{-4 \ell B D}}{8 \ell^2 B^2 D}
\frac{n}{\log n},\frac{n}{2B}\}$
the following holds: $M_k^B(m) \ge\alpha n$ implies
$M_k^B(m+\Delta) > \frac{\alpha}{2 B} e^{-2 \ell B D} n $ for every $0
\le\Delta\le D \frac{n}{k}$.
Then $\Pr(\cL(\alpha,B,D))\ge1-n^{-1}$.
\end{lemma}
\begin{pf}
As in the proof of Lemma~\ref{lemcreateGC}, it suffices to consider
fixed $m$ and~$k$,
and show that conditional on $\cF_m$, if $G(m)$ satisfies $M_k^B(m)\ge
\alpha n$,
then with probability at least $1-n^{-4}$ we have
$M_k^B(m+\Delta) > \frac{\alpha}{2B} e^{-2 \ell B D} n$
for every $0 \le\Delta\le\tilde{\Delta}$, where $\tilde{\Delta
}=\lfloor D n/k \rfloor$.

Condition on $\cF_m$, and let $C_1,\ldots,C_{r}$ be the components of $G(m)$
with sizes between $k$ and $Bk-1$.
Note that $r \ge M_k^B(m)/(Bk) \ge\alpha n/(Bk)$.

Starting from $G(m)$, we now analyze the next $\tilde{\Delta}$ steps.
We say that $C_{i}$ is \textit{safe} if in each of these steps none of
the $\ell$
randomly chosen vertices is contained in $C_{i}$, and we denote by $X$
the number of safe components.
Using $|C_{i}| \le Bk \le n/2$, note that $C_{i}$ is safe with probability
\[
(1-|C_{i}|/n)^{\ell\tilde{\Delta}} > e^{-2\ell\tilde{\Delta}|C_{i}|/n}
\ge e^{-2\ell BD} ,
\]
which gives $\E X \ge r e^{-2 \ell BD}$.
Clearly, the random variable $X$ can be written as $X=f(\underline
{v}_{m+1}, \ldots, \underline{v}_{m+\Delta})$,
where the $\underline{v}_{j}$ denote the $\ell$-tuples generated by the
$\ell$-vertex process in each step
(uniformly and independently).
The function~$f$ satisfies $|f(\omega)-f(\tilde{\omega})| \le\ell$ whenever
$\omega$ and $\tilde{\omega}$ differ in one coordinate.
So, using $r \ge\alpha n/(Bk)$, McDiarmid's inequality~\cite{McD}
implies that
$\Pr(X \le r e^{-2 \ell B D}/2)$ is at most
\[
\exp\biggl(-\frac{2[r e^{-2 \ell B D}/2]^2}{\tilde{\Delta} \ell
^2}\biggr)
\le\exp\biggl(-\frac{\alpha^2e^{-4 \ell B D}}{2\ell^2B^2D}\frac
{n}{k}\biggr)
\le n^{-4} .
\]

Suppose that $X > r e^{-2 \ell B D}/2$.
Since every safe component contributes at least $k$ vertices to every
$M_k^B(m+\Delta)$
with $0 \le\Delta\le\tilde{\Delta}$ (in each step all edges which
can be added are disjoint
from safe components), using $r \ge\alpha n/(Bk)$ we deduce that for
all such $\Delta$
we have
$M_k^B(m+\Delta) \ge kX > \alpha e^{-2 \ell B D} n/(2B)$, and the proof
is complete.
\end{pf}

Note that by considering instead the number $Y$ of vertices in safe
components one can prove the
slightly stronger bound $M_k^B(m+\Delta) > (1-\eps)\alpha e^{-2 \ell B
D} n $,
for $k$ not too large.

We are now ready to prove Theorem~\ref{th1}.
\begin{pf*}{Proof of Theorem~\ref{th1}}
Let $h_L(n)$ and $h_m(n)$ be nonnegative functions satisfying
$h_L(n)=o(n)$ and $h_m(n)=o(n)$,
and let $\delta>0$ be constant.
Let $\cX=\cX_n(\delta,h_L,h_m)$ denote the event that there exist
$m_1$ and $m_2$ satisfying $L_1(m_1) \le h_L(n)$,
$L_1(m_2) \ge\delta n$,
and $m_2 \le m_1+ h_m(n)$, so our aim is to show that $\Pr(\cX)\to0$
as $n\to\infty$.
We shall define a ``good'' event $\cG=\cG_n(\delta)$ such that $\Pr
(\cG
)\to1$ as $n\to\infty$
and show \textit{deterministically} that there is some $n_0$
such that for $n\ge n_0$, when $\cG$ holds, $\cX$ does not.

To be totally explicit, set $\alpha=\delta/4$, $A=5/\alpha^{\ell-1}$
and $D=1$.
Set $B=\ceil{2A\ell^2/\delta}$, and let
$\beta=\alpha e^{-2 \ell B}/(2B)>0$.
Finally, let $K=B^{1+\ceil{1/\beta}}$, noting that $K$ does not depend
on $n$.

Let $\cG$ be the event that $\cC(1)$, $\cC(\delta/4)$
and $\cL(\delta/4,B,D)$ all hold simultaneously.
By Lemmas~\ref{lemcreateGC} and~\ref{lemBinSize}, $\Pr(\cG)\ge
1-3n^{-1}=1-o(1)$.
The definition of $\cG$ ensures that if $n$ is large enough (larger
than some constant depending only on $\delta$ and $\ell$),
then for all $m\le5n$ and $k\le K$ the following hold:
\[
\mbox{(i) }N_{\ge k}(m)\ge\delta n/4 \quad\mbox{implies}\quad
\mbox{(ii) } L_1(m+\floor{An/k})\ge\delta n/(4\ell^2)
\]
and
\[
\mbox{(iii) }M_k^B(m)\ge\delta n/4 \quad\mbox{implies}\quad\mbox{(iv) }
M_k^B(m')\ge\beta n \quad\mbox{for all }m\le m'\le m+n/k .
\]

Suppose that $\cG$ holds, and that $m^-=\max\{m\dvtx L_1(m)\le h_L(n)\}$ and
$m^+=\min\{m\dvtx L_1(m)\ge\delta n\}$
differ by at most $h_m(n)$. It suffices to show deterministically that
if $n$ is large enough, then this leads to a contradiction.

Since $N_1(0)=n$ and $\cC(1)$ holds, we have $L_1(4n)\ge n/\ell^2$.
If $n$ is large enough, it follows that $m^-\le4n$, so $m^+\le5n$.

For\vspace*{1pt} $k\le K/B$ set $m_k=m^+-\delta n/(\ell^2 k)$, which is
easily seen to be positive; we ignore the irrelevant rounding to integers.
Since\vspace*{1pt} at most ${\ell\choose2}(m^+-m_k)<\ell^2(m^+-m_k)/2$ edges are
added passing
from $G(m_k)$ to $G(m^+)$, the components of $G(m_k)$ with size at most $k$
together contribute at most $k\ell^2(m^+-m_k)/2 \le\delta n/2$ vertices
to any one component of $G(m^+)$. It follows that
\[
N_{\ge k}(m_k)\ge L_1(m^+)-\delta n/2 \ge\delta n/2 .
\]
Suppose that $N_{\ge Bk}(m_k)\ge\delta n/4$. Then (i)
holds at step $m_k$ with $Bk\le K$ in place of $k$,
so (ii) tells us that by step
\[
m^*=m_k+\floor{An/(Bk)}\le m_k+\delta n/(2\ell^2 k) = m^+-\delta
n/(2\ell^2 k) = m^+-\Theta(n),
\]
we have $L_1(m^*) > \delta n/(4\ell^2)$, which is larger than $h_L(n)$
if $n$ is large enough.
Since $m^+-m^-\le h_m(n)=o(n)$, if $n$ is large enough we have $m^*<m^-$,
contradicting the definition of $m^-$.

It follows that $M_k^B(m_k)=N_{\ge k}(m_k)-N_{\ge Bk}(m_k)\ge\delta
n/4$. Using (iii) implies~(iv),
this gives $M_k^B(m^+)\ge\beta n$. Applying this for $k=1,B,B^2,\ldots
,B^{\ceil{1/\beta}}$
shows that $G(m^+)$ has more than $n$ vertices, a contradiction.
\end{pf*}

Setting\vspace*{1pt} $D=2\delta/\ell^2$ (instead of $D=1$), the proof
above shows that the number of steps between $m^-=\max\{m\dvtx
L_1(m)\le\delta/(4\ell^2)n\}$ and $m^+=\min\{ m\dvtx\break L_1(m)\ge\delta
n\}$ is at least $\delta n/(2\ell^2B^{\ceil{1/\beta}}) = f(\delta) n$,
where $f(\delta)$ essentially grows like the inverse of a double
exponential in $\delta^{-(\ell-1)}$ for $\delta\to0$.

\section{Results for merging rules}\label{secedge}

Although Theorem~\ref{th1} applies to any $\ell$-vertex rule, for many
questions, this class is too broad. Indeed, consider a~rule which only
joins two components when forced to (i.e., when presented with $\ell$
vertices from distinct components) and then joins the two smallest
components presented. Such a rule will \textit{never} join two of the
$\ell-1$ largest components, and it is not hard to see that during the
process $\ell-1$ giant components [with order $\Theta(n)$] will
emerge and grow simultaneously, with their sizes keeping roughly in
step. In what follows we could replace ``the largest
component'' by ``the union of the $\ell-1$ largest components'' and work
with arbitrary $\ell$-vertex rules, but this seems rather unnatural.

By an \textit{$r$-Achlioptas rule} we mean an $\ell$-vertex rule
with $\ell=2r$ that always joins (at least) one of the pairs
$\{v_{1},v_{2}\},\{v_{3},v_{4}\},\ldots,
\{v_{\ell-1},v_{\ell}\}$. (How we treat the case where one or more of these
pairs is in fact a single vertex will not be relevant.)
An \textit{Achlioptas rule} is an $r$-Achlioptas rule for any $r\ge1$.
Taking $r=2$ and insisting that only one edge is added gives the
original class of rules suggested by Achlioptas.

Let us say that an $\ell$-vertex rule is \textit{merging} if, whenever
$C$, $C'$ are distinct components with $|C|,|C'|\ge\eps n$,
then in the next step we have probability at least $\eps^\ell$ of joining
$C$ to $C'$. This\vspace*{-1pt} implies that the
probability that they are \textit{not} united after $m$ further steps
is at most $e^{-\eps^\ell m}$.
[We could replace $\eps^\ell$ by any $f(\eps)>0$,
and it suffices if the chance of merging in one of the next
few steps, rather than the next step, is not too small.]
Clearly, any Achlioptas rule is merging; with probability
at least $\eps^\ell$ all $r=\ell/2$ potential edges join $C$ to $C'$.
There are other interesting examples of merging rules (see Section~\ref{size}).

For merging rules we have the following variant of Lemma~\ref{lemcreateGC}.
We write $V_{\ge k}(m)$\vspace*{1pt} for the union of all components with size at
least $k$ in $G(m)$, so $|V_{\ge k}(m)|=N_{\ge k}(m)$.
%
%
\begin{lemma}\label{cGC2}
Let $\cR$ be a merging $\ell$-vertex rule, let $\eps>0$,
let $k\ge1$ and $m$ be integers, and set
$\Delta=2\ceil{\frac{2^\ell}{\eps^{\ell-1}}\frac{n}{k}}$.
Conditioned on $\cF_m$, with probability at least $1-\ell\exp(-c
n/k)$ there
is a component of $G(m+\Delta)$ containing at least $N_{\ge k}(m)-\eps
n$ vertices from $V_{\ge k}(m)$,
where $c=c(\eps,\ell)>0$.
\end{lemma}
\begin{pf}
Let $W=V_{\ge k}(m)$, so $|W|=N_{\ge k}(m)$.
We may assume that $|W|-\eps n \ge0$. Let $\alpha=|W|/n \ge\eps$.
Until the point that there are $\ell-1$ components between them
containing at least
$(\alpha-\eps/2)n$\vspace*{1pt} vertices from $W$, at each step we have probability
at least $\alpha(\eps/2)^{\ell-1}$ of choosing $\ell$ vertices of $W$
in distinct components to form $\bv_j$, in which case the number of components
meeting~$W$ must decrease by (at least) one. As in the proof of
Lemma~\ref{lemcreateGC}, it
follows that off an event whose probability
is exponentially small in $n/k$, after $\Delta/2$ steps we do have
$\ell-1$ components $C_1,\ldots,C_{\ell-1}$ together containing at
least $(\alpha-\eps/2)n$
vertices of $W$. Ignoring any containing fewer than $\eps n/(2\ell)$ vertices
of $W$, using the property of merging rules noted above, the probability
that some pair of the remaining $C_i$ are not joined in the next
$\Delta/2$
steps is exponentially small in $n/k$.
\end{pf}

It is easy to check that we may take $c(\eps,\ell)=\eps/\ell^\ell$.
With this technical result in hand, we now prove Theorem~\ref{thgu}.
\begin{pf*}{Proof of Theorem~\ref{thgu}}
We outline the argument, much of which is very similar to the proof of
Theorem~\ref{th1} given in the previous section.

Let $\eps> 0$ be given and set $\delta=\eps/5$. Lemma~\ref{cGC2} implies
that there is some $A=A(\delta,\ell)$ such that for any fixed $k$, it
is very likely
that (i) there is a~component of $G(m+\lfloor An/k\rfloor)$ containing
at least
$N_{\ge k}(m)- \delta n$\vspace*{1pt} vertices. 
By Lemma~\ref{lemBinSize}, for every fixed $B$ there is some $\beta
=\beta(\delta,\ell,B)>0$ such that
if (ii)\vspace*{1pt} $M_k^B(m)=N_{\ge k}(m)-N_{\ge Bk}(m)\ge\delta n$,\vspace*{1pt} then it is
very likely
that (iii) $M_k^B(m')\ge\beta n$ for all $m\le m'\le m+n/k$,
say.\vadjust{\goodbreak}

To be more precise, let $B=\ceil{A\ell^2/\delta}$
and $K=B^{1+\ceil{1/\beta}}$. Then it follows
easily from Lemma~\ref{lemBinSize}, Lemma~\ref{cGC2} and the union bound
that for $n$ large enough there is a good event $\cG=\cG_n(\delta)$ such
that $\Pr(\cG)\to1$ and such that whenever~$\cG$ holds,
then for all $m\le n^2$ and $k\le K$, (i) holds and (ii) implies (iii).

Suppose that $\cG$ holds and that $m^+=\min\{m\dvtx N_{\ge K}(m) \ge
L_1(m)+\eps n\}$ exists.
It suffices to show deterministically that if $n$ is large enough, then
this leads to a contradiction.
Since $\cG$ holds, considering (i) with $m=0$ and $k=1$ shows that for
some $C=C(\delta,\ell)$
we have $L_1(Cn)\ge(1-\delta)n > (1-\eps) n$, so $m^+\le Cn$.

For $k\le K/B$, set $m_k=m^+ - 2\delta n/(\ell^2 k)$.
Recall that $V_{\ge k}(m)$ denotes the the union of all components with
size at least $k$ in $G(m)$.
Since at most ${\ell\choose2}(m^+-m_k) < \delta n/k$ edges are added
passing from
$G(m_k)$ to $G(m^+)$, vertices outside of $V_{\ge k}(m_k)$ contribute
at most $2\delta n$ vertices to $V_{\ge k}(m^+)$. Hence,
\[
N_{\ge k}(m_k) \ge N_{\ge k}(m^+) - 2\delta n \ge N_{\ge K}(m^+) -
2\delta n \ge L_1(m^+)+(\eps-2\delta)n .
\]
Suppose that $N_{\ge Bk}(m_k) \ge N_{\ge k}(m_k)- \delta n$.
Then (i) (with $Bk$ in place of $k$) tells us that by step
\[
m= m_k+\floor{An/(Bk)} \le m_k+\delta n/(\ell^2 k) = m^+-\delta
n/(\ell
^2 k) < m^+
\]
there exists a component of $G(m)$ containing at least
\[
N_{\ge Bk}(m_k)-\delta n \ge N_{\ge k}(m_k)- 2\delta n \ge
L_1(m^+)+(\eps-4\delta)n
> L_1(m^+)
\]
vertices, which contradicts $G(m^+) \supseteq G(m)$.
It follows that $M_k^B(m_k)\ge\delta n$. Using (ii) implies (iii)
we deduce that $M_k^B(m^+)\ge\beta n$. Applying this for
$k=1,B$, $B^2,\ldots,B^{\ceil{1/\beta}}$
and counting vertices in $G(m^+)$ gives a contradiction.
\end{pf*}

Working\vspace*{1pt} through the conditions on the constants in the proof above,
and using $D=3\delta/\ell^2$ instead of $D=1$ when applying
Lemma~\ref{lemBinSize}, one can check
that for some positive constants $c$ and $d$ depending only on $\ell$
the result
holds for any $\eps=\eps(n)\ge d/(\log\log n)^{1/(\ell-1)}$,
with $K=K(\eps)\le\exp(\exp(c \eps^{-(\ell-1)}))$.

%
\begin{theorem}\label{thcont}
Let $\cR$ be a merging $\ell$-vertex rule.
For each $n$, let\break $(G(m))_{m\ge0}=(G_n^{\cR}(m))_{m\ge0}$ be the random
sequence of graphs on $\{1,2,\ldots,n\}$ associated to $\cR$.
Given any function $h_m(n)$ that is $o(n)$,
and any constants $0\le a<b$,
the probability that there exist $m_1$ and $m_2$
with $L_1(G(m_1))\le an$, $L_1(G(m_2))\ge bn$ and $m_2\le m_1+h_m(n)$
tends to $0$ as $n\to\infty$.
\end{theorem}

Note that for merging rules, Theorem~\ref{thcont} implies
the conclusion of Theorem~\ref{th1}; a~``jump'' from $o(n)$ to $\ge
\delta n$
implies a ``jump'' from $\le\delta n/2$ to $\ge\delta n$.
\begin{pf*}{Proof of Theorem~\ref{thcont}}
Let $a<b$ be given, and set $\eps=(b-a)/2$. Using Theorem~\ref{thgu}
we
may assume that
there exists $K=K(\eps, \ell)$ such that $N_{\ge K}(m) < L_1(m)+\eps n$
for all $m$.\vadjust{\goodbreak}
Suppose that $m^-=\max\{m\dvtx L_1(m)\le an\}$ and $m^+=\min\{m\dvtx L_1(m)\ge
bn\}$ differ by at
most $h_m(n)$. Set $m^*=m^+ - \eps n/(2\ell^2 K)$. As before, we have
%
%
\begin{equation}\label{eqthcont}
N_{\ge K}(m^*)\ge L_1(m^+)- \ell^2 K (m^+-m^*) > (b - \eps) n =
(a+\eps
)n .
\end{equation}
If $n$ is large enough, which we assume, then $m^+ \le m^- + h_m(n)$
implies $m^*<m^-$.
This gives $N_{\ge K}(m^*)
\le N_{\ge K}(m^-)
< L_1(m^-)+\eps n \le(a + \eps)n$, contradicting (\ref{eqthcont}).
\end{pf*}

Let us remark that Theorem~\ref{thcont} (which can be proved without
first proving Theorem~\ref{thgu}) gives an alternative proof of
Spencer's ``no two giants'' conjecture; if at any time there are two
components with at least $\eps n$ vertices, then in the step after the
last such time, $L_1$ must increase by at least $\eps n$ in a~single
step. Hence, Theorem~\ref{thcont} implies that if $\cR$ is merging,
then for any $\eps>0$ we have $\max_m L_2(m)\le\eps n$ whp.
%
%
\begin{corollary}\label{corcont}
Let $\cR$ be a merging $\ell$-vertex rule.
If $\cR$ is globally convergent, then $\rho^{\cR}$ is continuous on
$[0,\infty)$.
\end{corollary}
\begin{pf}
Let $\rho(t)=\rho^{\cR}(t)$. We have $0\le\rho(t)\le{\ell\choose2}t$,
so $\rho$ is continuous at $0$. Suppose $\rho$ is discontinuous at
some $t>0$.
Since $\rho$ is increasing, $\sup_{t'<t}\rho(t')<\inf_{t'>t}\rho(t')$,
so we may pick $a<b$ with $\sup_{t'<t}\rho(t')<a<b<\inf_{t'>t}\rho(t')$.
By definition of global convergence,
for any fixed $\eps>0$,
%
%
\begin{equation}\label{b1}
\Pr\bigl( L_1\bigl(\floor{(t-\eps)n}\bigr)\le an \mbox{ and } L_1\bigl(\floor
{(t+\eps
)n}\bigr)\ge bn
\bigr)\ge1-\eps,
\end{equation}
if $n$ is large enough. It follows
as usual that there is some $\eps(n)\to0$ such that~(\ref{b1}) holds
with $\eps=\eps(n)$. But this contradicts Theorem~\ref{thcont}.
\end{pf}

\section{Convergence considerations}\label{secconv}

From the beginning, a key question about Achlioptas processes
has been which rules are globally convergent. In some cases,
local convergence has been established, but as far as we are aware,
global convergence has not been shown for any nontrivial rules.

We now turn to the proof of Theorem~\ref{thlg},
that local convergence implies global convergence
for merging rules (in particular, for Achlioptas rules). We comment
further on local convergence below. Theorem~\ref{thlg} is easy
to deduce from Theorem~\ref{thgu}; we shall give a more direct
proof that seems more informative.
\begin{pf*}{Proof of Theorem~\ref{thlg}}
Suppose $\cR$ is locally convergent. Then there exist
functions $\rho_k\dvtx [0,\infty)\to[0,1]$ such that (\ref{ktconv}) holds
for any fixed $k \geq1$ and $t \geq0$.
Since $N_k$ changes by at most $2k$ when an edge is added to a~graph,
it follows easily that each $\rho_k$ is continuous (indeed Lipschitz).
From monotonicity of the underlying process, it is easy to see that
for each $k$, the function $\rho_{\le k}(t)=\sum_{j\le k}\rho_j(t)$
is decreasing.

Define $\rho=\rho^{\cR}$ by
\[
\rho(t) = 1-\sum_{k=1}^\infty\rho_k(t) = 1-\lim_{k\to\infty}
\rho_{\le
k}(t) ,
\]
so $\rho\dvtx[0,\infty)\to[0,1]$ is increasing.
We claim that for any fixed $t>0$ and $\eps>0$, the probability that
%
%
\begin{equation}\label{claim}
\sup_{0\le t'<t} \rho(t')-\eps\le\frac{L_1(\floor{tn})}{n} \le
\rho
(t)+\eps
\end{equation}
tends to $1$ as $n\to\infty$. This clearly implies that $L_1(\floor
{tn})/n\pto\rho(t)$
whenever~$\rho$ is continuous at $t$, which is the definition of global
convergence.
Corollary~\ref{corcont} then implies that $\rho$ is continuous.

The upper bound in (\ref{claim}) is immediate; by definition of $\rho$
there is some $K$ such that $\rho_{\le K}(t) \ge1-\rho(t)-\eps/4$.
Summing (\ref{ktconv}) up to $K$ gives $N_{\le K}(\floor{tn})/n \ge
1-\rho(t)-\eps/2$
whp. When $n$ is large enough, this bound implies $L_1(\floor
{tn})/n\le
\rho(t)+\eps$.

For the lower bound, we combine the ``sprinkling'' argument
of Erd\H os and R\'enyi~\cite{ERgiant} with Lemma~\ref{cGC2}. Choose $t'<t$
such that $\rho(t')$ is within $\eps/2$ of the supremum,
and let $m_1=\floor{t'n}$ and $m_2=\floor{tn}$, so $m_2-m_1=\Theta(n)$.
It suffices to show that $L_1(m_2)/n\ge\rho(t')-\eps/2$ holds whp.
In doing so we may assume that \mbox{$\rho(t')-\eps/2\ge0$}.
For any constant $K$, whp we have $N_{\le K}(m_1)/n\le\rho_{\le
K}(t')+\eps/4 \le1-\rho(t')+\eps/4$,
so $N_{\ge K}(m_1)/n\ge\rho(t')-\eps/4$ whp. If $K$ is large enough
(depending only
on $t'$ and $\eps$), Lemma~\ref{cGC2}
then gives $L_1(m_2)/n \ge\rho(t')-\eps/2$ whp, as required.
\end{pf*}
%
%
\begin{remark}
Since nonmerging $\ell$-vertex rules have received some attention
(see, e.g.,~\cite{NLT2011}), let us spell out what our method gives for
such rules. Lem\-ma~\ref{cGC2} applies in this case provided ``there is a
component containing'' is changed to ``there are $\ell-1$ components
together containing.''
Let $L(m)$ denote the sum of the sizes of the $\ell-1$ largest components
in $G(m)$. With this modified Lemma~\ref{cGC2}, the proof of
Theorem~\ref{thgu}
goes through with $L_1$ replaced by $L$. The same is true of
Theorem~\ref{thcont}
[with an extra $-(\ell-2)K$ in (\ref{eqthcont}),
since the largest $\ell-1$ components may not
all be large]. Finally, Corollary~\ref{corcont} and Theorem~\ref{thlg}
similarly go through, now with $\rho$ defined using $L$ rather than $L_1$.
\end{remark}

\section{Size rules}\label{size}

So far, even in the Achlioptas-rule case our rules have been very general,
making choices between the given edges using any information about the
current graph. There is a natural much smaller class (of vertex or
Achlioptas rules)
called \textit{size rules}, where only the sequence $c_1,\ldots,c_\ell$
of the orders of the components containing the presented vertices
$v_1,\ldots,v_\ell$ may be used to decide which edge(s) to add. (Here
we suppress the dependence on the step $m$ in the notation.)
Note that the product rule is a size rule.

In fact, most past results concern \textit{bounded size rules}; here
there is a constant $B$ such that all sizes $c_i> B$ are treated the
same way by the rule, so the rule only ``sees'' the data
$(\min\{c_i,B+1\})_{i=1}^{\ell}$. Perhaps the simplest example is the
Bohman--Frieze process, the bounded size rule with $B=1$ in which the
edge $v_1v_2$ is added if $c_1=c_2=1$, and otherwise $v_3v_4$ is added.
Bohman and Frieze~\cite{BF2001} showed that for a closely related rule
there is no giant component when $m\sim0.535n$. [The actual rule they
used considered whether $v_1$ and $v_2$ are isolated in the graph
formed by all pairs \textit{offered} to the rule, rather than the graph
$G(m)$ formed by the pairs \textit{accepted} so far.]

Considering, for simplicity, rules in which one edge is added
at each step, a~key property of bounded size rules is that at each step,
the expected change in $N_k$ can be expressed as a simple function
of $N_1,N_2,\ldots,N_{\max\{k,B\}}$. (It is clear that the rate
of formation of $k$-vertex components can be so expressed;
for the rate of destruction, consider separately
the cases $k$ joins to $k'$ for each $k'\le B$ and the
case $k$ joins to some $k'>B$.)
Spencer and Wormald~\cite{SpencerWormald2007}, who
considered bounded size Achlioptas rules,
and Bohman and Kravitz~\cite{BohmanKravitz2006},
who considered a large subset of such
rules, noted that
in this case one can easily use Wormald's ``differential equation
method''~\cite{WormaldDEM}
to show that the rule is locally convergent,
and that the $\rho_k(t)$ satisfy certain differential equations.
This remark applies to all bounded size $\ell$-vertex rules.

Resolving a conjecture of Spencer~\cite{potp},
Spencer and Wormald~\cite{SpencerWormald2007} proved that
any bounded size 2-Achlioptas rule exhibits a phase
transition: there is some $\tc$, depending on the rule,
such that for $t<\tc$, whp $L_1(\floor{tn})=o(n)$ [in fact $O(\log n)$],
while for $t>\tc$, $L_1(\floor{tn})=\Omega(n)$ whp.
They conjectured that any bounded size 2-Achlioptas rule is globally convergent,
and that the phase transition is second order (continuous).
Theorem~\ref{thlg}
establishes both these conjectures.

Very recently, Janson and Spencer~\cite{JansonSpencer2010} established
bounds on the size of the giant component in the Bohman--Frieze process
just above the (known) critical point $\tc$. They deduce that
\textit{if} it is globally convergent, then the right derivative of
$\rho$ at $\tc$ has a certain specific value. The required ``if'' part
is established by Theorem~\ref{thlg}.

Informally, let us call a size rule \textit{nice} if there is some $K$
such that, for each~$k$, the expected change in $N_k$ is a function
of $N_1,N_2,\ldots,N_{\max\{k,K\}}$. [More precisely,
the individual decisions whether to create or destroy a component
of size $k$ depend only on the data $(\min\{c_i,k'+1\})_{i=1}^\ell$ where
$k'=\max\{k,K\}$ and $c_i$ is the size of the component containing
$v_i$.] Just as in the bounded size
case, using the differential equation method,
it is easy to show that any nice rule
is locally convergent. Hence, by Theorem~\ref{thlg}, any nice
\textit{merging} rule
is globally convergent with continuous phase transition;
this applies to all \textit{nice} Achlioptas rules.

The simplest examples of nice rules have $K=1$, that is, only compare
component sizes.
One example is ``join the two smallest.''
For $\ell=3$ this rule is mentioned briefly by
Friedman and Landsberg~\cite{FriedmanLandsberg2009} as another
example of a rule that should be explosive,
and discussed by D'Souza and Mitzenmacher~\cite{DSouzaMitzenmacher2010},
who ``established'' the explosive nature of the transition for this
and a related nice rule numerically; Theorem~\ref{th1}
contradicts these predictions.

Another nice rule is the following: join the smaller of $C_1$ and $C_2$
to the smaller of $C_3$ and $C_4$, where $C_i$ is the component
containing $v_i$. We call this the ``dCDGM'' rule since it was
introduced by da Costa, Dorogovtsev, Goltsev and Mendes~\cite{dCDGM}.
Note that this is \textit{not} an Achlioptas rule, but it \textit{is}
merging; if $|C|,|C'|\ge\eps n$ then with probability at least $\eps^4$
we choose $v_1,v_2\in C$ and $v_3,v_4\in C'$ and so join $C$ to $C'$.
Hence, the dCDGM rule, which is locally convergent by the differential
equation method, is globally convergent and has a continuous phase
transition. Da Costa, Dorogovtsev, Goltsev and Mendes~\cite{dCDGM}
proposed this rule as simpler to analyze than the product rule, but at
least as likely to have a discontinuous phase transition. For a brief
discussion of their arguments, see the end of the
\hyperref[secintro]{Introduction}.

There are many open questions concerning the precise nature
of the phase transitions in various Achlioptas and related processes.
One of the most intriguing is the following:
Is the product rule globally convergent?


%

\printaddresses

\end{document}